\documentclass[11pt]{amsart}

\usepackage{amssymb}
\usepackage{amsmath}
\usepackage{amsthm}
\usepackage{tikz}

\newtheorem{thm}{Theorem}[section]

\newtheorem{corr}{Corollary}

\numberwithin{equation}{section}

\newcommand{\bmt}[1]{\mbox{\boldmath $#1$}}
\newcommand{\bg}{\bmt{\gamma}}
\newcommand{\bG}{\bmt{\Gamma}}
\newcommand{\bJ}{\bmt{J}}
\newcommand{\be}{\bmt{\beta}}
\newcommand{\bT}{\bmt{T}}
\newcommand{\bN}{\bmt{N}}
\newcommand{\bB}{\bmt{B}}
\newcommand{\bTt}{\tilde{\bmt{T}}}
\newcommand{\bNt}{\tilde{\bmt{N}}}
\newcommand{\bBt}{\tilde{\bmt{B}}}
\newcommand{\man}{\mathcal{M}}
\newcommand{\re}{\mathbb{R}}
\newcommand{\sph}{\mathbb{S}}
\newcommand{\ex}{\textnormal{exp}}

\begin{document}

\title[Integrable geodesic flows on tubular sub-manifolds]{Integrable geodesic flows on tubular sub-manifolds}

\author{Thomas Waters}
\address{University of Portsmouth}
\email{thomas.waters@port.ac.uk}

\begin{abstract}In this paper we construct a new class of surfaces whose geodesic flow is integrable (in the sense of Liouville). We do so by generalizing the notion of tubes about curves to 3-dimensional manifolds, and using Jacobi fields we derive conditions under which the metric of the generalized tubular sub-manifold admits an ignorable coordinate. Some examples are given, demonstrating that these special surfaces can be quite elaborate and varied.
\end{abstract}


\maketitle

\section{Introduction}

Let $(\man^n,g)$ be a Riemannian manifold. We may write the geodesic equations as a Hamiltonian system on the cotangent bundle $T^*\man$. The integral curves of this vector field constitute a flow, called the geodesic flow. If this Hamiltonian system has $n$ first integrals, functionally independent and in involution, we say the geodesic flow is Liouville integrable (see, for example, \cite{kling},\cite{josesaletan},\cite{sakai},\cite{docarmoriemm}). Throughout this paper we will describe this property with phrases such as ``the manifold has integable geodesic flow'', or ``the manifold is integrable''.

The known surfaces and manifolds with integrable geodesic flow is surprisingly small. The classic examples are surfaces of revolution and ellipsoids (see \cite{kling} for a detailed treatment or \cite{tabach} for an unusual approach), both known for many years (18th and 19th centuries respectively). The case of ellipsoids has been extended to quadratic manifolds in general \cite{tabach}, and there was a recent series of new examples of integrable manifolds in the form of Lie groups (see \cite{bolsinov},\cite{bols2} for reviews and references). On the other hand, in some recent papers \cite{TWspherical},\cite{TWmonge},\cite{TWetds} the author has shown that there are classes of surfaces and manifolds whose geodesic flow is {\it not} integrable. It would seem that surfaces/manifolds with integrable geodesic flow are indeed very rare and special, and often those that are known are not of a type which can be visualized easily \cite{dullin}. In this paper we will add to the list of integrable surfaces, by first considering the geodesic flow on tubes about curves and then generalizing this notion. We will show that in the examples we construct the Hamiltonian has an ignorable coordinate and hence by Noether's theorem \cite{josesaletan} admits a linear integral; together with the Hamiltonian itself this means the surface has integrable geodesic flow.

We begin with an observation about curves in $\re^3$. Let $\bg:I\subset\re\to\re^3$ be a simple smooth regular curve, and let $\{\bmt{T},\bmt{N},\bmt{B}\}$ be the associated Frenet frame. The tube of radius $\rho_0$ about $\bg$ has parameterisation \[ \bmt{\beta}:I\times\mathbb{S}^1\to\re^3:\bmt{\beta}(s,\psi;\rho_0)=\bg(s)+\rho_0\cos\psi\bmt{N}(s)+\rho_0\sin\psi\bmt{B}(s).\] The line element is found to be \begin{equation} d\varsigma^2=\big[(1-k_1(s)\rho_0\cos\psi)^2+k_2(s)^2\rho_0^2\big]ds^2+2k_2(s)\rho_0^2 dsd\psi+\rho_0^2d\psi^2, \label{ler3} \end{equation} where $k_1(s),k_2(s)$ are the curvature scalars (curvature and torsion) associated with $\bg$, and we restrict $\rho_0$ such that the tube is a well-defined surface (see later). From the line element we can find the Hamiltonian for the geodesic flow, $\mathcal{H}=\tfrac{1}{2}g^{ij}p_ip_j$, and we observe that $s$ only enters the Hamiltonian explicitly via the curvature scalars $k_{1,2}$. Thus if $k_{1,2}$ are constants, $s$ will be an ignorable coordinate leading to a linear integral, $p_s$, and the geodesic flow will be integrable. Hence the geodesic flow on the tube about a helix in $\re^3$ is integrable.

Let us state this explicitly: the tube around a curve in $\re^3$ has integrable geodesic flow if the curve has constant curvature scalars. There are a number of ways of generalizing this statement which spring to mind, and some are discussed in the conclusions. In this work we will focus on the following: suppose the curve is instead embedded in some manifold $\man$ of dimension 3, and the curve has constant curvature scalars with respect to $\man$; is the geodesic flow on the tube about this curve integrable? We will see that the answer is yes when $\man$ is a space form, but also in the case that $\man$ is of a more general class and {\it not} of constant curvature.

First some background theory which will be useful in what follows. The Frenet-Serret equations are generalized to curves in Riemannian manifolds in the following way: let $\bg=\bg(s):I\subset\re\to\man$ be a simple smooth regular curve parameterized by arc-length. Letting $\bmt{T}=\dot{\bg}$ and $\{\bmt{T},\bmt{N},\bmt{B}\}$ be an orthonormal frame in $T_{\bg} \man$, the frame evolves along $\bg$ according to (\cite{Gutkin},\cite{Barros},\cite{Tamura}): \[ \frac{D}{ds}\begin{pmatrix}
  \bmt{T}\\ \bmt{N}\\ \bmt{B}
\end{pmatrix}=\begin{pmatrix}
  0 & k_1(s) & 0 \\ -k_1(s) & 0 & k_2(s) \\ 0 & -k_2(s) & 0
\end{pmatrix}\begin{pmatrix}
  \bmt{T}\\ \bmt{N}\\ \bmt{B}
\end{pmatrix}. \] Here and throughout we will let $D/ds$ denote the covariant derivative with respect to $\dot{\bg}$, as opposed to $D$ in \cite{Gutkin}, $\bar{\nabla}_{\bmt{T}}$ in \cite{Barros} and $D_{\dot{\bg}}$ in \cite{Tamura}. $k_1$ is often known as the `geodesic curvature', however we will simply refer to $k_{1,2}$ as the `curvature scalars'. If $k_{1,2}$ are constants we will refer to $\bg$ as a `constant curvature curve' (as opposed to `helix', `proper helix' etc. \cite{Barros},\cite{Tamura}).

It will be convenient to denote by $T_{\bg}^{\perp}\man$ the orthogonal complement of $\dot{\bg}$ in $T_{\bg}\man$ (i.e.\ span$(\bmt{N},\bmt{B})$). In Euclidean space, the tube about a curve $\bg$ is the locus of circles of fixed radius orthogonal to $\bg$. However, in $\man$ the tube about $\bg$ is the locus of {\it geodesic} circles whose radial geodesics have tangent vectors in $T_{\bg}^{\perp}\man$ at $\bg$. As such a parameterisation of the tube of radius $\rho_0$ about $\bg\in\man$ would be \[ \bmt{\beta}:I\times\mathbb{S}^1\to\man:\bmt{\beta}(s,\psi;\rho_0)=\textrm{exp}_{\bg(s)}\big(\rho_0\cos\psi\bmt{N}(s)+\rho_0\sin\psi\bmt{B}(s)\big), \] where $\ex_p:T_p\man\to\man$ denotes the exponential map defined as follows: let $\bG$ be the unique geodesic through $p$ with tangent vector $\bmt{v}$ at $p$; then $\ex_p(\bmt{v})=\bG(\|\bmt{v}\|)$. As alluded to previously, we restrict $\bg$ and $\rho_0$ as required in order for the tube to be well defined, by avoiding either the curvature of $\bg$ being too high at a certain point, reaching conjugate/focal points along radial geodesics (see \cite{sakai},\cite{TWbif}), or $\bg$ passing close to itself leading to self-intersections.

The question we wish to address is: for what $\man$ does the tube $\bmt{\beta}(s,\psi;\rho_0)$ have integrable geodesic flow? While tubes as sub-manifolds have been much studied previously (see for example \cite{gray}), it is typically in the context of the area and volume of the tube, and not in terms of the geodesic flow {\it on} the tube. In the following section we formulate the problem using Jacobi fields and prove that if $\man$ is a space form then the geodesic flow on the tube about a constant curvature curve is integrable; we also give some examples. An advantage to formulating the problem in terms of Jacobi fields is we lay the foundations for further study of evolutes, focal surfaces, etc. (see \cite{TWbif} for example).  In Section 3 we consider curves in more general manifolds not of constant curvature, and find some conditions for the geodesic flow on tubular sub-manifolds to be integrable; again we give some examples. In each case considered we show the existence of an ignorable coordinate, i.e. the metric admits a Killing vector field and the sub-manifold is therefore invariant under a 1-parameter group of isometries; integrability follows naturally. Section 4 contains some conclusions and further comments.

\section{Curves in space forms}

Consider the following parameterisation of the neighbourhood of $\bg$: \[ \be:I\times\mathbb{S}^1\times \re^+\to\man:\be(s,\psi,\rho)=\textnormal{exp}_{\bg(s)}\big(\rho\cos\psi\bmt{N}+\rho\sin\psi\bmt{B}\big) \] with $\be(\rho=0)=\bg$ and the appropriate restriction on $\rho$ as discussed previously. We could describe this as a `geodesic cylindrical coordinate system', with $\bg$ as the axis and $\rho=\rho_0$ defining the tube of radius $\rho_0$ about $\bg$. We will show that under certain conditions the line element on $\rho=\rho_0$ is independent of $s$; this will imply the geodesic flow on $\rho=\rho_0$ is integrable. For this we require the pairwise inner products of $d\bmt{\beta}/ds$ and $d\bmt{\beta}/d\psi$.

\begin{figure}
\scalebox{0.8}{
\begin{tikzpicture}
\clip (2,0.6) rectangle (14,8);
\draw[thick] (1,0) .. controls (9,3) and (1,6) .. (3,8);
\draw (4,4.8) arc (130:111:7);
\draw (4.38,4.15) arc (120:100:7);
\draw (4.6,3.5) arc (110:91:7);
\draw[thick] (8,0) .. controls (16,3) and (8,6) .. (10,8);
\draw (5,2) node {$\bg(s)$};
\draw (4.1,3.3) node {\footnotesize{$s_{-1}$}};
\draw (4,4) node {\footnotesize{$s_0$}};
\draw (3.7,4.6) node {\footnotesize{$s_1$}};
\draw (8,5) node {$\bG(\rho;s_0,\psi_0)$};
\draw[->,thick] (4.38,4.15) -- (3.5,6);
\draw[->,thick] (6.65,4.98) -- (6,6.85);
\draw (4.5,6) node {$\bmt{T}(s_0)$};
\draw (7.6,6.8) node {$\bmt{J}_s(\rho;s_0,\psi_0)$};
\draw (11.38,4.15) arc (160:135:2.5);
\draw (11.38,4.15) arc (128:106.5:6);
\draw (11.38,4.15) arc (90:72.5:5);
\draw[dashed,rotate around={20:(11.38,4.15)}] (11,3.3) arc (-100:100:2.4cm and 0.8cm);
\draw[->, thick] (13.4,5.175) -- (11.5,5.9);
\draw (12.7,6.5) node {$\bmt{J}_\psi(\rho;s_0,\psi_0)$};
\draw (13.2,4.8) node {\footnotesize{$\psi_0$}};
\draw (11.8,5.3) node {\footnotesize{$\psi_1$}};
\draw (13.2,3.7) node {\footnotesize{$\psi_{-1}$}};
\draw[->,thick] (4.26,2) -- (4.325,2.1);
\draw[->,thick] (11.26,2) -- (11.325,2.1);
\draw[->,thick] (2.825,6.6) -- (2.795,6.65);
\draw[->,thick] (9.825,6.6) -- (9.795,6.65);
\draw[dashed] (7,3) .. controls (6.88,5) and (6.7,5) .. (5.4,7);
\draw[->] (7.5,4.7) arc (-30:-150:0.77);
\end{tikzpicture}
}
\caption{The radial geodesics along $\bg$ give two variations through geodesics by letting either $s$ or $\psi$ vary; the variation fields $\bJ_s$ and $\bJ_{\psi}$ are shown.}\label{varfields}
\end{figure}
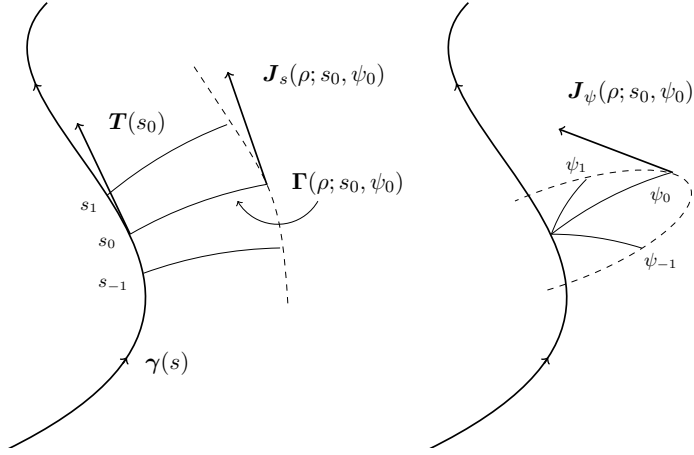

 If we fix $s=s_0,\psi=\psi_0$ and consider the radial geodesic \[ \bG(\rho;s_0,\psi_0)=\ex_{\bg(s_0)}\big(\rho\cos\psi_0\bmt{N}(s_0)+\rho\sin\psi_0\bmt{B}(s_0)\big), \] we may describe two variation fields along this geodesic (see Figure \ref{varfields}) \[ \bJ_s(\rho)=\frac{\partial \bmt{\beta}}{\partial s},\quad\textrm{and}\quad\bJ_{\psi}(\rho)=\frac{\partial \bmt{\beta}}{\partial \psi} \]
 where we will suppress cumbersome notation such as $\bJ_s(\rho;s_0,\psi_o)$. Both $\bJ_s$ and $\bJ_{\psi}$ solve the Jacobi equation (here and throughout a prime denotes $d/d\rho$)  \begin{equation} \frac{D^2}{d\rho^2}\bJ_*=R(\bJ_*,\bG')\bG',\qquad (*=s,\psi), \label{jaceq}\end{equation} as indeed do all Jacobi fields; what distinguishes one Jacobi field from the next is the initial conditions. Focussing on $\bJ_s$ first, we note \[ \bJ_s(\rho=0)=\frac{\partial \bmt{\beta}}{\partial s}(\rho=0)=\frac{\partial}{\partial s}\big( \bmt{\beta}(\rho=0)\big)=\frac{\partial}{\partial s}\big(\bg\big)=\bmt{T}. \] Also (letting $|_0$ denote $\rho=0$), \[ \left.\frac{D\bJ_s}{d\rho}\right|_0=\left.\left[\frac{D}{d\rho}\left(\frac{\partial \bmt{\beta}}{\partial s}\right)\right]\right|_0=\left.\left[\frac{D}{ds}\left(\frac{\partial \bmt{\beta}}{\partial \rho}\right)\right]\right|_0=\frac{D}{ds}\left[\left.\left(\frac{\partial \bmt{\beta}}{\partial \rho}\right)\right|_0\right]. \] Remembering that $\bmt{\beta}(s,\psi,\rho)$ denotes the unit speed geodesic, parameterized by $\rho$, passing through $\bg(s)$ with tangent vector $\cos\psi\bmt{N}(s)+\sin\psi\bmt{B}(s)$, we have \[ \frac{D}{ds}\left[\left.\left(\frac{\partial \bmt{\beta}}{\partial \rho}\right)\right|_0\right]=\frac{D}{ds}\left(\cos\psi\bmt{N}+\sin\psi\bmt{B}\right)=\cos\psi(-k_1\bmt{T}+k_2\bmt{B})+\sin\psi(-k_2\bmt{N}).  \] In summary both $\bJ_s$ and $\bJ_\psi$ solve the Jacobi equation \eqref{jaceq} with initial conditions \begin{equation} \bJ_s(0)=\bmt{T},\quad \frac{D\bJ_s}{d\rho}(0)=-k_1\cos\psi\bmt{T}-k_2\sin\psi\bmt{N}+k_2\cos\psi\bmt{B}, \label{jsdata} \end{equation} and \begin{equation} \bJ_\psi(0)=\bmt{0},\quad \frac{D\bJ_\psi}{d\rho}(0)=-\sin\psi\bmt{N}+\cos\psi\bmt{B}. \label{jpsidata}  \end{equation} We are now ready to prove the first theorem.

\begin{thm}
  Let $\bg=\bg(s):I\subset \re\to\man$ be a simple smooth regular curve in the 3-dimensional Riemannian manifold $\man$. If $\man$ is a space form $(\re^3/\sph^3/\mathbb{H}^3)$ and $\bg$ is a constant curvature curve, then the geodesic flow on the tube about $\bg$ is integrable.
\end{thm}

\begin{proof}
  If $\man$ is a space form with constant sectional curvature $K_0$, $\bG$ is a unit speed geodesic of $\man$, and $\bJ_*$ is a Jacobi field orthogonal to $\bG'$, then \cite{docarmoriemm}  \[ R(\bJ_*,\bG')\bG'=-K_0\bJ_* \] and hence the Jacobi equation \eqref{jaceq} becomes \[ \frac{D^2\bJ_*}{d\rho^2}=-K_0\bJ_*. \] Given w.l.o.g.\ $K_0=-1,0,1$ we may simply solve this Jacobi equation for $\bJ_s$ and $\bJ_\psi$ with the data given in \eqref{jsdata},\eqref{jpsidata}, then take their pairwise inner products to find the line element on the tube $\rho=\rho_0$, which is \[ d\varsigma^2=\big[(F_0-k_1(s)G_0\cos\psi)^2+k_2(s)^2G_0^2\big]ds^2+2k_2(s)G_0^2 ds d\psi+G_0^2 d\psi^2 \] where \[ \left. \begin{array}{c}
    F_0=\cos(\rho_0),\ G_0=\sin(\rho_0) \\
    F_0=1,\ G_0=\rho_0 \\
    F_0=\cosh(\rho_0),\ G_0=\sinh(\rho_0)
  \end{array} \right\} \quad \textrm{if}\quad \left\{\begin{array}
    {c} K_0=1,\\ K_0=0, \\ K_0=-1.
  \end{array} \right. \]
   In each case the coordinate $s$ only enters the line element via the curvature scalars $k_1$ and $k_2$, and hence if $\bg$ is a constant curvature curve $s$ will be an ignorable coordinate and the geodesic flow on the tube will be integrable.
\end{proof}

Before giving an example we can extend this result by observing that a `tube', defined as an equidistant surface, consists of a circle in $T_{\bg}^{\perp}\man$ projected under exp into $\man$; however we can also consider the `generalized tube' where any simple closed curve, fixed w.r.t.\ $s$, in $T_{\bg}^{\perp}\man$ is projected into $\man$. That is, if $(f(\psi),g(\psi))$ is a simple closed curve in the plane, then the generalized tube (intuitively understood as this simple closed curve being ``carried along'' by $\bg$) also has integrable geodesic flow.

\begin{corr}
  The geodesic flow on the generalized tube \[ \be(s,\psi,\rho_0)=\textnormal{exp}_{\bg(s)}\big(\rho_0 f(\psi)\bmt{N}(s)+\rho_0 g(\psi)\bmt{B}(s) \big) \] has integrable geodesic flow if $\bg\subset\man$ has constant curvature scalars and $\man$ is $\re^3,\sph^3$ or $\mathbb{H}^3$.
\end{corr}

We give an example. Using the Hopf parameterisation of $\sph^3$, \[ \bmt{\sigma}(\eta,\theta,\phi)=(\sin\eta\cos\theta,\sin\eta\sin\theta,\cos\eta\cos\phi,\cos\eta\sin\phi), \] we consider the curve $(\eta,\theta,\phi)=(\eta_0,\alpha t,\beta t)$. We require $\alpha/\beta\in\mathbb{Q}$ (in order to avoid the tube self-intersecting for $t\in\mathbb{R}$) in which case $\bg$ is a knot on the Clifford torus in $\sph^3$. $t$ is only the arc-length if $\sqrt{\alpha^2\sin^2\eta_0+\beta^2\cos^2\eta_0}\equiv p=1$. The curvature scalars are found to be (correcting \cite{Tamura}) \[ k_1=\frac{(\alpha^2-\beta^2)\sin\eta_0\cos\eta_0}{p^2}, \quad k_2=\frac{\alpha \beta}{p^2}. \] Note this curve has constant curvature scalars, and hence its tube has integrable geodesic flow.

\begin{figure}[h!]
\begin{center}
 \includegraphics[width=0.75\textwidth]{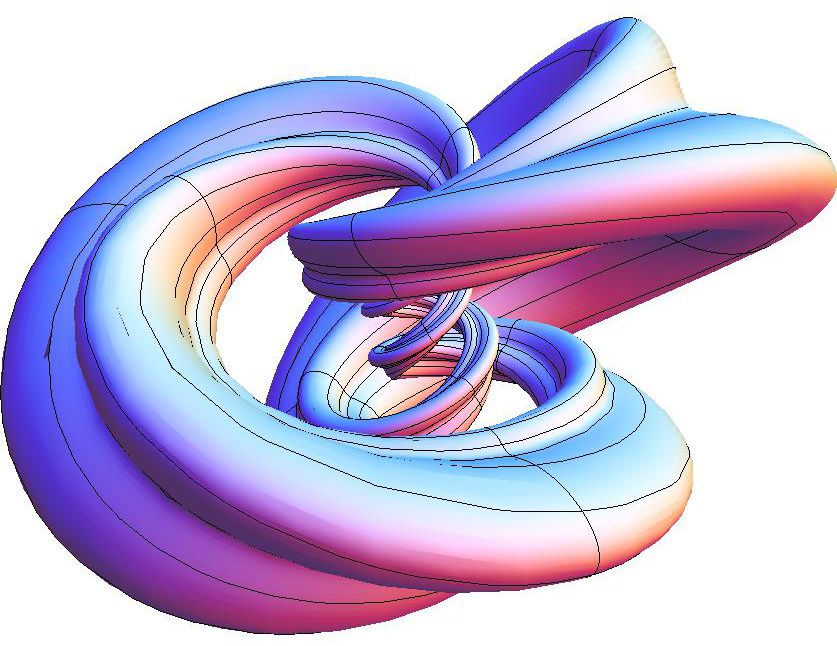}\caption{The stereograhic projection from $\sph^3$ into $\re^3$ of a generalized tube about a constant curvature curve (see text for details). This surface has integrable geodesic flow.} \label{stereo}
\end{center}
\end{figure}

In Fig \ref{stereo} we show the stereographic projection from $\sph^3$ into $\re^3$ of the generalized tube with $f=(1+0.3\cos(3\psi))\cos\psi,\ g=(1+0.3\cos(3\psi))\sin\psi,\rho_0=0.2$ about the constant curvature curve described in the proceeding paragraph with $\alpha=5,\beta=2,\eta_0=\pi/4$.

\section{Curves in more general manifolds}

We suspect that there may be manifolds more complicated than the space forms which admit integrable sub-manifolds in the form of (generalized) tubes about curves if the curves take advantage of some symmetry in the manifold. Indeed this is the case as we will show.

\begin{thm}
  Let $\bg=\bg(s):I\subset \re\to\man$ be a simple smooth regular curve in the 3-dimensional Riemannian manifold $\man$. The neighbourhood of $\bg\subset\man$ can be parameterized as follows: \[ \be:I\times\mathbb{S}^1\times \re^+\to\man:\be(s,\psi,\rho)=\textnormal{exp}_{\bg(s)}\big(\rho\cos\psi\bmt{N}+\rho\sin\psi\bmt{B}\big) \] with $\be(\rho=0)=\bg$ and the appropriate restriction on $\rho$. If in the tubular neighbourhood of $\bg$ the sectional curvatures are independent of $s$, and $\bg$ is a constant curvature curve, then the tube $\be(s,\psi;\rho_0)$ will have integrable geodesic flow.
\end{thm}

\begin{proof}
  We begin by extending the $\bT,\bN,\bB$ frame into the neighbourhood of $\bg$ by parallel transport along the radial geodesics $\bG$, i.e.\ define $(\bTt,\bNt,\bBt)(s,\psi,\rho)$ via \[ \nabla_{\bG'}(\bTt,\bNt,\bBt)\equiv\frac{D}{d\rho}(\bTt,\bNt,\bBt)=\bmt{0},\quad (\bTt,\bNt,\bBt)(\rho=0)=(\bT,\bN,\bB). \] Parallel transport preserves lengths and angles so $(\bTt,\bNt,\bBt)$ is an orthonormal basis in the neighbourhood of $\bg$, and we can write \[ \bJ_*=t_*\bTt+n_*\bNt+b_*\bBt \] where everything depends on $s,\psi,\rho$, and the Jacobi equation separates into three ODE's: \begin{align*}
    t_*''+\langle R(\bG',\bJ_*)\bG',\bTt\rangle=0, \\ n_*''+\langle R(\bG',\bJ_*)\bG',\bNt\rangle=0, \\ b_*''+\langle R(\bG',\bJ_*)\bG',\bBt\rangle=0.
  \end{align*} Note $\bG'(s,\psi,\rho)=\cos\psi\bNt+\sin\psi\bBt$ (since $\nabla_{\bG'}(\cos\psi\bNt+\sin\psi\bBt)=\bmt{0}$ and $\bG'(s,\psi,0)=\cos\psi\bN+\sin\psi\bB$), thus each of these equations has 12 terms. Many are zero however (using the symmetries of the Riemann tensor and \cite{sakai}), reducing to \begin{align}
    t_*''+t_*\cos^2\psi(\bNt,\bTt,\bNt,\bTt)+t_*\sin^2\psi(\bBt,\bTt,\bBt,\bTt)=0,\nonumber \\
    n_*''+(\bNt,\bBt,\bNt,\bBt)[n_*\sin^2\psi-b_*\cos\psi\sin\psi]=0, \\
    b_*''+(\bNt,\bBt,\bNt,\bBt)[b_*\cos^2\psi-n_*\cos\psi\sin\psi]=0=0,\nonumber
  \end{align} (using the notation for sectional curvature  in \cite{docarmoriemm}) with initial data \begin{align*} t_s(0)=1,\ n_s(0)=0,\ b_s(0)=0,\\ t_s'(0)=-k_1(s)\cos\psi,\ n_s'(0)=-k_2(s)\sin\psi,\ b_s'(0)=k_2(s)\cos\psi, \end{align*} and \begin{align*} t_{\psi}(0)=0,\ n_{\psi}(0)=0,\ b_{\psi}(0)=0,\\ t_{\psi}'(0)=0,\ n_{\psi}'(0)=-\sin\psi,\ b_{\psi}'(0)=k_2(s)\cos\psi. \end{align*} In general, since the coefficients of these differential equations depend on $s,\psi,\rho$ then the solutions will depend on $s,\psi,\rho$. However, if the sectional curvatures are independent of $s$, then the coefficients in these differential equations do not depend explicitly on $s$. Hence the dependence of $t_*$ etc.\ on $s$ only manifests itself via the initial conditions, which in turn only depend on $s$ via the curvature scalars $k_{1,2}$. Thus if these scalar functions are in fact constants, then $t_*$ etc.\ will not depend explicitly on $s$ and therefore $\langle\bJ_s,\bJ_s\rangle=t_s^2+n_s^2+b_s^2$ etc.\ will not depend on $s$. Hence $s$ will be an ignorable coordinate for the Hamiltonian function describing the geodesic flow on the tubes $\be(s,\psi;\rho_0)$, which will therefore be integrable.
\end{proof}

\noindent We give an example. Consider the manifold (a degenerate 3-ellipsoid) \[ \frac{x_1^2}{a^2}+\frac{x_2^2}{a^2}+\frac{x_3^2}{b^2}+\frac{x_4^2}{b^2}=1. \] Modifying the Hopf coordinates in the obvious way, the metric tensor is \[ diag\big(a^2\cos^2\eta+b^2\sin^2\eta,a^2\sin^2\eta,b^2\cos^2\eta \big) \] and the curve $(\eta,\theta,\phi)=(\eta_0,\alpha t,\beta t)$ has constant curvatures \[ k_1=\left(\frac{b^2\beta^2}{p^2}-\frac{a^2\alpha^2}{p^2}\right)\ \frac{\cos\eta_0\sin\eta_0}{q},\quad k_2=\frac{a b \alpha \beta}{p^2 q} \] and Frenet frame \[ \bT=(0,\alpha/p,\beta/p),\quad \bN=(1/q,0,0),\quad \bB=\left(0,\frac{b \beta\cos\eta_0}{a p},-\frac{a \alpha \tan\eta_0}{b p}\right), \] where $p=\sqrt{a^2\alpha^2\sin^2\eta_0+b^2\beta^2\cos^2\eta_0}$ and $q=\sqrt{a^2\cos^2\eta_0+b^2\sin^2\eta_0}$.

Now to show that the tube around $(\eta,\theta,\phi)=(\eta_0,\alpha s/p,\beta s/p)$ is integrable we need to show that the sectional curvatures in the neighbourhood of $\bg$ do not depend on $s$. The difficulty is that Theorem 2 is phrased in terms of the geodesic cylindrical coordinates $s,\psi,\rho$, whereas our manifold is parameterised by $\eta,\theta,\phi$ and the connection between these two coordinate systems is not known explicitly. However, we do know that if we let $R$ denote a generic component of the Riemann tensor, then \begin{equation} \frac{\partial R}{\partial s}=\frac{\partial R}{\partial \eta}\frac{\partial \eta}{\partial s}+\frac{\partial R}{\partial \theta}\frac{\partial \theta}{\partial s}+\frac{\partial R}{\partial \phi}\frac{\partial \phi}{\partial s}. \label{genR} \end{equation} The last two terms are zero as the metric tensor depends only on $\eta$. To show the first term is zero we consider $\partial \eta/\partial s$.

Suppose we take a point $s=s_0$ on $\bg$ and consider the radial geodesic whose initial tangent vector makes the angle $\psi_0$ in the $\bN(s_0),\bB(s_0)$ plane; this geodesic will solve the initial value problem \[ \ddot{\eta}+\Gamma^1_{11}\dot{\eta}^2+\Gamma^1_{22}\dot{\theta}^2+\Gamma^1_{33}\dot{\phi}^2=0,\quad \ddot{\theta}+2\Gamma^2_{12}\dot{\eta}\dot{\theta}=0,\quad \ddot{\phi}+2\Gamma^3_{13}\dot{\eta}\dot{\phi}=0, \] with \begin{align*} \big(\eta,\theta,\phi,\dot{\eta},\dot{\theta},\dot{\phi}\big)(0)=\left(\eta_0,\frac{\alpha s_0}{p},\frac{\beta s_0}{p},\frac{\cos\psi_0}{q},\frac{\sin\psi_0 b \beta\cos\eta_0}{a p},-\frac{\sin\psi_0 a \alpha \tan\eta_0}{b p}\right). \end{align*} However, since the Christoffel symbols do not depend on $\theta,\phi$, we see the following $4$-$d$ first order system decouples: \[ \dot{\eta}=u,\ \dot{u}=-\Gamma^1_{11}u^2-\Gamma^1_{22}v^2-\Gamma^1_{33}w^2,\ \dot{v}=-2\Gamma^2_{12}uv,\ \dot{w}=-2\Gamma^3_{13}uw, \] with \[ \eta(0)=\eta_0,\ u(0)=\frac{\cos\psi_0}{q},\ v(0)=\frac{\sin\psi_0 b \beta\cos\eta_0}{a p},\ w(0)=-\frac{\sin\psi_0 a \alpha \tan\eta_0}{b p}. \] This system is completely independent of $s$, as moving along $\bg$ only changes $\theta(0),\phi(0)$ which are not part of this decoupled system. As such, the solutions for $(\eta,u,v,w)$ are independent of $s$.

Since $d\eta/ds=0,$ this means via \eqref{genR} that the elements of the Riemann tensor are independent of $s$ and hence the sectional curvatures in the neighbourhood of $\bg$ are independent of $s$. Hence Theorem 2 implies the tube around this curve is integrable.

\section{Conclusion}

We have shown that there is a broad class of sub-manifolds (which could be described as `generalized tubes') which have geodesic flow that is integrable in the sense of Liouville. This adds to the very limited number of integrable surfaces known to date, and we can see from Figure \ref{stereo} that these sub-manifolds can be quite varied and elaborate. What's more our formulation in terms of Jacobi fields provides a foundation for further study of caustics and envelopes in curved manifolds. Our results arose out of generalizing to $\man$ the following observation in $\re^3$: the tube around a curve in $\re^3$ is integrable if the curve has constant curvature scalars. As mentioned in the Introduction there are other possibilities in generalizing this result.

Firstly we could consider changing the previous statement to: the tube around a curve in $\re^3$ is integrable {\it if and only if} the curve has constant curvature scalars. An obvious place to begin would be to consider the tube around an ellipse. Rigourous proofs of non-integrability using differential Galois theory could be attempted as in \cite{TWspherical},\cite{TWmonge},\cite{TWetds}, however strong evidence of non-integrability can be provided much more quickly by the method of Poincar\'{e} section (see Appendix A). While it seems very likely the ``if and only if'' statement is true, it would be very hard to prove: there is no topological obstruction to integrability (the tubes are either cylinders or torii) and as such the geometry of the tube would need to be specified exactly (or at best within a large family of curves with parameters).

A second generalization would be to consider curves in $\re^n$, and construct the following manifold: at each point on the curve take the orthogonal complement to the tangent vector and place an $(n-2)$-dimensional sphere. This would generate an $(n-1)$-dimensional manifold, and we would suspect that this manifold might be integrable if the $n-1$ curvature scalars associated with the curve are constants. This is however too much to ask: we find that again $s$ is an ignorable coordinate leading to an integral of the motion, which together with the ``energy'' gives two integrals, but for curves in $\re^4$ or above this would not be enough to imply integrability. Again, we offer numerical evidence: for the 3-dimensional manifold around a constant curvature curve in $\re^4$, we use the integral $p_s$ to reduce the geodesic flow to a 2-dimensional Hamiltonian system and take a Poincar\'{e} section. The results are very like those of Appendix A and will not be shown. Nonetheless there is some numerical evidence to suggest that there may be families of integrable manifolds in this broad class and this may be worth further investigation.

\appendix

\section{Poincar\'{e} section of the geodesic flow on the tube about an ellipse in $\re^3$}

The method of Poincar\'{e} section is quite standard in Hamiltonian dynamics (\cite{josesaletan},\cite{TWspherical}) and we only sketch the method here, in the context of a tube about an ellipse. Suppose the phase space of a Hamiltonian system has coordinates $(s,\psi,p_s,p_\psi)$. Fixing the numerical value of the Hamiltonian (for geodesic flows we may choose an arc-length parameterisation), we choose initial conditions on some hyperplane of phase space (for example $s=0$) and follow a trajectory recording subsequent intersections with the hyperplane in the same direction (i.e.\ $s=0\ \textrm{mod}\ 2\pi,p_s>0)$. This will produce a series of points in the $(\psi,p_\psi)$-plane. According to the Liouville-Arnol'd theorem, if the Hamiltonian system is integrable these points will trace out a closed curve, and if not they will fill a region of the plane.

\begin{figure}[h]
\begin{center}
 {\includegraphics[width=0.5\textwidth]{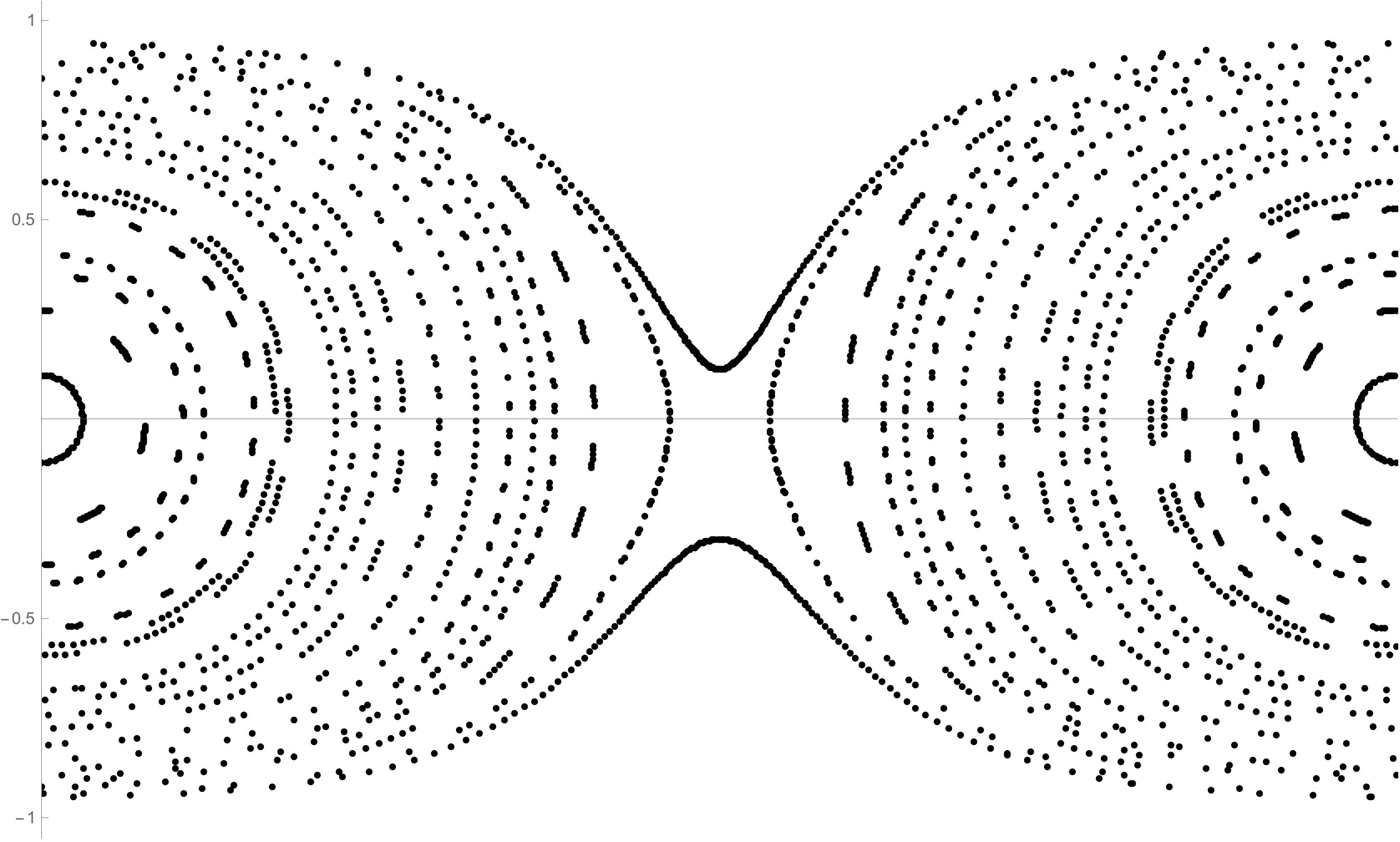}\includegraphics[width=0.5\textwidth]{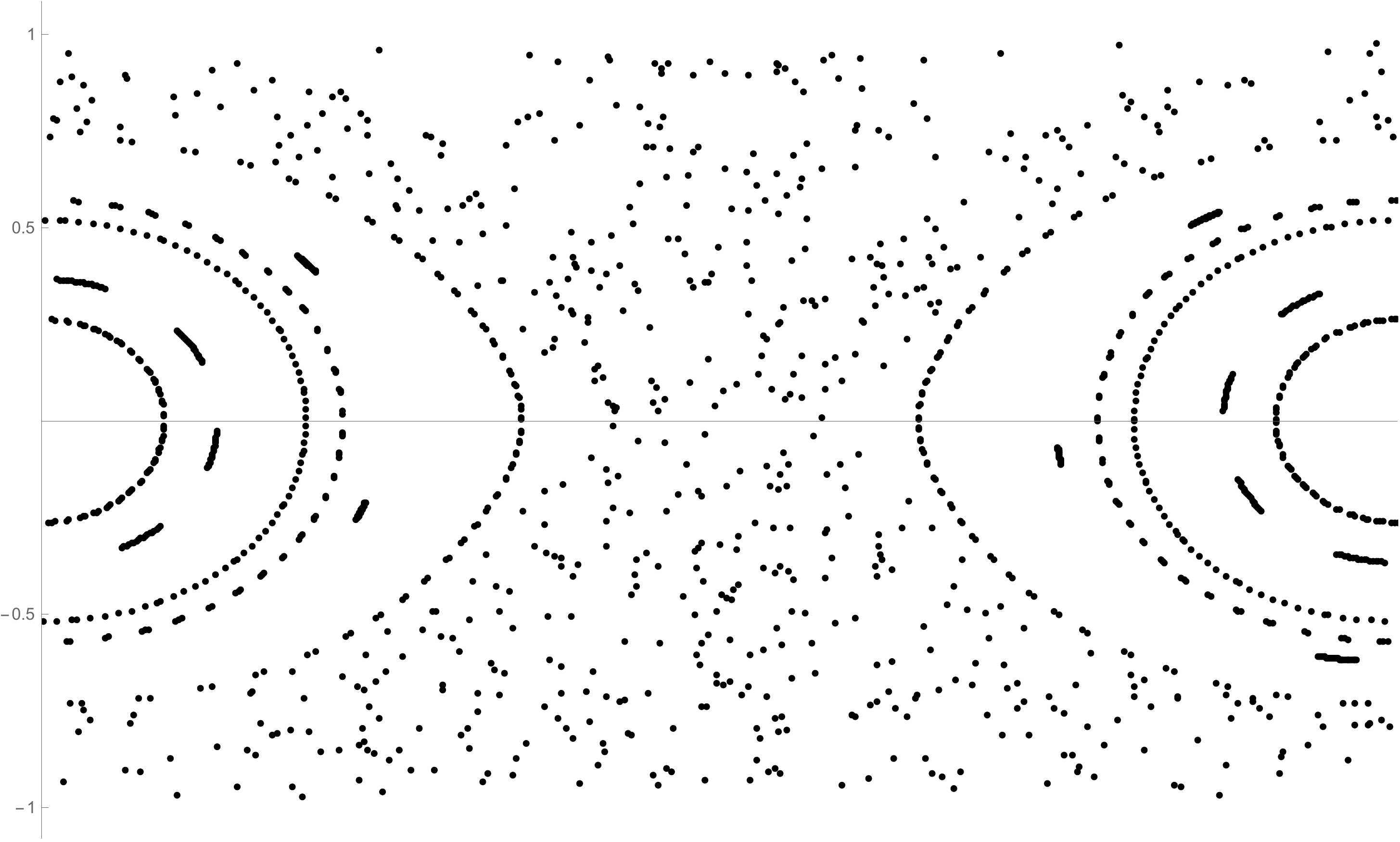}}
 \caption{Two sets of Poincar\'{e} sections on the tube about ellipses (see text for details). The horizontal axes are $\psi\in(0,2\pi)$, and the vertical $p_\psi\in(-1,1)$.} \label{psecs}
\end{center}
\end{figure}

In Figure \ref{psecs} we show two sets of Poincar\'{e} sections for a variety of initial conditions. On the left is the tube of radius 1 about the ellipse with semi-axes $(2,2)$ (i.e.\ a torus, which as we know has integrable geodesic flow), and on the right the tube of radius 1 about the ellipse with semi-axes $(2,2.5)$. We can easily see the disruption of the separatrix and the transition to chaotic orbits via the homoclinic tangle; clear evidence of non-integrability.

\bibliographystyle{proc_igc_plain}

\bibliography{tubebib}

\end{document}